\documentclass{article}
\usepackage[utf8]{inputenc}
\usepackage{amssymb, amsmath}
\usepackage{amsthm}
\usepackage{amsfonts} 
\usepackage{tikz-cd}
\title{Counting Semistable Representations of Quivers over Finite Fields}
\author{Jiuzhao Hua} 
\newtheorem{thm}{Theorem}[section]
\newtheorem{lem}{Lemma}[section]
\newtheorem{cor}{Corollary}[section]
\newtheorem{dfn}{Definition}[section]
\newtheorem{cjc}{Conjecture}[section]

\usepackage[parfill]{parskip}

\begin{document}
\date{\vspace{-0.5cm}}\date{} 
% Toggle commenting to test
\maketitle
\begin{abstract}
In this paper, we derive a closed formula for the number of isomorphism classes of absolutely indecomposable semistable representations of an arbitrary quiver 
over a finite field with a fixed dimension vector. This generalises a formula for Kac polynomials given by Hua. A key step in the proof is to show that any 
representation of a quiver with a nilpotent endomorphism over an arbitrary field admits a structured filtration by subrepresentations compatible with the nilpotent action.
\end{abstract}

\section{Introduction}

Let $\mathbb{N}$ denote the set of non-negative integers, $\mathbb{Z}$ the ring of integers, $\mathbb{Q}$ the field of rational numbers, 
$\mathbb{F}_q$ a finite field with $q$ elements, where $q$ is a prime power, and $\bar{\mathbb{F}}_q$ the algebraic closure of $\mathbb{F}_q$.

Given a positive integer $n$ and an $n \times n$ matrix $C = [a_{ij}]$ with non-negative integer entries, we define a \textit{quiver} $\Gamma$ associated with $C$ as follows:
\begin{itemize}
	\item The set of vertices is $\Gamma_0 =\{ 1,2,\ldots,n \}$.
	\item For each pair of vertices $(i, j)$, the quiver has $a_{ij}$ arrows from $i$ to $j$. The set of arrows is denoted by $\Gamma_1$.
\end{itemize}
We call $C$ the \textit{companion matrix} of the quiver $\Gamma$.

Let $k$ be an arbitrary field. A \textit{representation} \( V \) of the quiver \( \Gamma \) over \( k \) consists of a collection of \( k \)-vector spaces \( (V_i) \) for each vertex \( i \in \Gamma_0 \), 
and a collection of \( k \)-linear maps \( (V_a: V_i \to V_j) \) for each arrow \( a\!: i \to j \) in \( \Gamma_1 \).  
The vector \( (\dim V_1, \dots, \dim V_n) \in \mathbb{N}^n \) is called the \textit{dimension vector} of \( V \), denoted by \( \dim V \).  

A \textit{homomorphism} from a representation \( V = (V_i, V_a)_{i \in \Gamma_0,\/ a \in \Gamma_1} \) to another representation \( W = (W_i, W_a)_{i \in \Gamma_0,\/ a \in \Gamma_1} \)
is a collection of $k$-linear maps $(f_i: V_i \to W_i)_{i \in \Gamma_0}$ such that for each arrow $a\!: i \to j$ in $\Gamma_1$, the following diagram commutes:
\[
\begin{tikzcd}
V_i \arrow[r, "V_a"] \arrow[d, "f_i"'] & V_j \arrow[d, "f_j"] \\
W_i \arrow[r, "W_a"] & W_j
\end{tikzcd}.
\]

Let $\mathrm{Mat}_{s \times t}(k)$ denote the space of all $s \times t$ matrices over $k$, for $s,t \in \mathbb{N}$.
Let \( \mathrm{Rep}(\alpha, k) \) denote the space of all representations of \( \Gamma \) over \( k \) with dimension vector \( \alpha \in \mathbb{N}^n \).  
Then, we have
\begin{equation} \label{rep space}
\mathrm{Rep}(\alpha, k) = \bigoplus_{(i\to j) \in \Gamma_1} \mathrm{Hom} ( k^{\alpha_i}, k^{\alpha_j} ) \cong \bigoplus_{(i\to j) \in \Gamma_1} \mathrm{Mat}_{\alpha_j\times\alpha_i}(k).
\end{equation}

For \( m \in \mathbb{N} \), let \( \mathrm{GL}(m, k) \) be the general linear group of degree \( m \) over \( k \), and define  
\[
\mathrm{GL}(\alpha, k) := \prod_{i \in \Gamma_0} \mathrm{GL}(\alpha_i, k).
\]
This group acts naturally on \( \mathrm{Rep}(\alpha, k) \) as follows. Let \( V = (V_i, V_a)_{i \in \Gamma_0, a \in \Gamma_1} \) be a representation with
dimension vector \( \alpha \), and let \( g = (g_i)_{i \in \Gamma_0}$ be an element of  $\mathrm{GL}(\alpha, k) \). The action of \( g \) on \( V \) yields a new representation 
\( V' = (V'_i, V'_a)_{i \in \Gamma_0,\/ a \in \Gamma_1} \) defined by:
\begin{itemize}
    \item \( V'_i = V_i \) for each vertex \( i \in \Gamma_0 \).
    \item For each arrow \( a\!: i \to j \) in \( \Gamma_1 \), the map \( V'_a: V'_i \to V'_j \) is given by
    \[
        V'_a = g_j V_a g_i^{-1}.
    \]
\end{itemize}
As such, two representations are \textit{isomorphic} if and only if they lie in the same orbit under this group action.

For a comprehensive treatment of quiver representations, we refer the reader to \textit{An Introduction to Quiver Representations} by Harm Derksen and Jerzy Weyman~\cite{D-W 2017}.

Semistable representations were first introduced by King~\cite{AK 1994} in a geometric framework. In his foundational work, King used geometric invariant theory to define stability conditions 
for quiver representations. Subsequently, Rudakov~\cite{AR 1997} introduced stability conditions for abelian categories. Reineke~\cite{M-R 2003} later developed stability and slope functions 
for the category of quiver representations in a purely algebraic setting. This algebraic reformulation proved particularly fruitful, as it enabled combinatorial approaches to the study of these representations.

The number of isomorphism classes of absolutely indecomposable representations of a quiver over a finite field with a fixed dimension vector is encoded by the Kac polynomial. 
Hua \cite{JH 2000} obtained a closed formula for this polynomial (Theorem~4.6). In this paper, we establish an analogous closed formula for semistable representations.

Recall from Reineke \cite{M-R 2003}, each stability function $\theta\!: \mathbb{Z}^n \to \mathbb{Z}$ induces to a slope function $\mu$ defined by:
\[
\mu(\alpha) := \frac{\theta(\alpha)} {\operatorname{ht}(\alpha)}, \enspace\text{ where } \enspace \operatorname{ht}(\alpha) := \sum_{i=1}^n\alpha_i \text{ for } \alpha=(\alpha_1,\ldots,\alpha_n) \in \mathbb{Z}^n\backslash\{0\}.
\]
The \textit{slope} of a representation $M$ is defined as $\mu(M) := \mu(\dim M)$. Following Reineke \cite{M-R 2003},
a representation $M$ is called \textit{semistable} (respectively, \textit{stable}) if for every non-zero proper subrepresentation $U$ of $M$, we have $\mu(U) \le \mu(M)$ (respectively, $\mu(U) < \mu(M)$).

Fix a rational number $\mu \in \mathbb{Q}$. The full subcategory of all semistable representations of $\Gamma$ over $k$ of slope $\mu$ 
forms an abelian category that is closed under kernels, images, cokernels, and extensions. 

Let $S$ be a constructible algebraic set defined over $\mathbb{F}_q$, that is, $S$ is the subset of $(\bar{\mathbb{F}}_q)^n$ defined by a finite collection of equations and inequalities of the form
\begin{align*}
&f_i(x_1,\dots,x_n) = 0, \qquad 1 \le i \le s, \\
&g_j(x_1,\dots,x_n) \ne 0, \qquad 1 \le j \le t,
\end{align*}
where $f_i, g_j \in \mathbb{F}_q[x_1,\dots,x_n]$.  
We say that $S$ \emph{exhibits polynomial count behavior} if there exists a polynomial $p(t)\in\mathbb{Q}[t]$ such that, for all integers $k \ge 1$,
\[
\bigl|\{(x_1,\dots,x_n)\in S \mid x_i \in \mathbb{F}_{q^k}\ \text{for all } 1 \le i \le n\}\bigr| = p(q^k).
\]

Using the notation introduced above, it follows from \eqref{rep space} that
\begin{equation}
	\left| \mathrm{Rep}(\alpha, \mathbb{F}_q) \right| = q^{\sum_{1 \leq i, j \leq n} a_{ij} \alpha_i \alpha_j}.
\end{equation}

The order of the general linear group $\mathrm{GL}(m, \mathbb{F}_q)$ is given by
\begin{equation}\label{order of GL}
|\mathrm{GL}(m, \mathbb{F}_q)| = \prod_{i=0}^{m-1} (q^m - q^i).
\end{equation}
By convention, we set $\left| \mathrm{GL}(0, \mathbb{F}_q) \right| = 1$. Hence, the order of the product group $\mathrm{GL}(\alpha, \mathbb{F}_q)$ is
\begin{equation}
\left| \mathrm{GL}(\alpha, \mathbb{F}_q) \right| = \prod_{i=1}^n |\mathrm{GL}(\alpha_i, \mathbb{F}_q)| = \prod_{i=1}^n \prod_{s=0}^{\alpha_i - 1} (q^{\alpha_i} - q^s).
\end{equation}
It is well-known that the representation space $\mathrm{Rep}(\alpha, \mathbb{F}_q)$ and the linear group $\mathrm{GL}(\alpha, \mathbb{F}_q)$ 
exhibit polynomial count behavior. In what follows, we will treat their cardinalities as polynomials in $q$, without further mention of this property.

Recall that the Euler form associated with $\Gamma$ is defined by:
\begin{equation}\label{euler_form}
\langle \alpha, \beta \rangle := \sum_{i=1}^n\alpha_i\beta_i - \sum_{i,j=1}^na_{ij}\alpha_i\beta_j \text{ for } \alpha, \beta \in \mathbb{Z}^n.
\end{equation}

Fix a rational number $\mu\in \mathbb{Q}$, and let 
\[
	\Delta_\mu =\{\/ \alpha\in\mathbb{Z}^n \mid \theta(\alpha) = \mu \cdot \mathrm{ht}(\alpha) \/\}.
\]
Then $\Delta_{\mu}$ is a free abelian group. Let $\Delta_\mu^+ = \Delta_\mu \cap \mathbb{N}^n$. 
For $\alpha\in\Delta_\mu^+$, define
\[
	\mathrm{R}_\mu^{ss}(\alpha, \mathbb{F}_q) = \{\/ M \in \mathrm{Rep}(\alpha, \mathbb{F}_q) \mid M \text{ is semistable} \/\}.
\]
It is clear that every nonzero representation in $\mathrm{R}_\mu^{ss}(\alpha, \mathbb{F}_q)$ has slope~$\mu$. 
A dimension vector $d \in \mathbb{N}^n$ is said to be \emph{coprime with respect to the slope function $\mu$} if
\[
\mu(e) \neq \mu(d) \quad \text{for all dimension vectors } 0 < e < d.
\]
If $\alpha$ is coprime with respect to $\mu$, then Reineke’s
Theorem~3.8~\cite{M-R 2003} guarantees the existence of stable representations of
dimension vector $\alpha$. In particular,
\[
\mathrm{R}_\mu^{ss}(\alpha, \mathbb{F}_q) \neq \varnothing .
\]

Reineke \cite{M-R 2003} shows that $\mathrm{R}_\mu^{ss}(\alpha, \mathbb{F}_q)$ exhibits polynomial count 
behavior, and its cardinality is given by the following formula:
\begin{equation}\label{num_ss}
\frac{|\mathrm{R}_\mu^{ss}(\alpha, \mathbb{F}_q)|} {|\mathrm{GL}(\alpha, \mathbb{F}_q)|} = \sum_{\alpha_*}(-1)^{s-1} q^{-\sum_{i<j}\langle\alpha^j, \alpha^i\rangle}\prod_{i=1}^s
\frac{|\mathrm{Rep}(\alpha^i, \mathbb{F}_q)|} {|\mathrm{GL}(\alpha^i, \mathbb{F}_q)|} ,
\end{equation}
where the sum runs over all $s$-tuples $\alpha_*=(\alpha^1,\ldots,\alpha^s)$ of non-zero dimension vectors with $s\ge 1$, $\sum_{i=1}^s\alpha^i = \alpha$, and 
$\mu(\sum_{i=1}^k \alpha^i) > \mu(\alpha)$ for $1\le k < s$. It follows from formula \ref{num_ss} that $|\mathrm{R}_\mu^{ss}(\alpha, \mathbb{F}_q)| \in \mathbb{Z}[q]$.

For $\alpha\in\Delta_\mu^+$, let $A_\mu^{ss}(\alpha, q)$ (respectively, $A_\mu^{s}(\alpha, q)$) denote the number of isomorphism classes of absolutely indecomposable semistable (respectively, stable) 
representations of $\Gamma$ over $\mathbb{F}_q$ with dimension vector $\alpha$. 
By convention, $A_\mu^{ss}(\alpha, q) = A_\mu^{s}(\alpha, q) = 0$ if $\alpha = 0$.

It was shown by Reineke~\cite{M-R 2006} that $A_\mu^{s}(\alpha,q)$ is a polynomial in $q$ with integer coefficients; that is, there exists a polynomial $f(t)\in\mathbb{Z}[t]$ such that
\[
A_\mu^{s}(\alpha,q)=f(q)
\]
for all prime powers $q$. 
Mozgovoy and Reineke~\cite{M-R 2014} conjectured that $A_\mu^{ss}(\alpha,q)$ is a polynomial in $q$ with non-negative integer coefficients.
As shown in Section~3, we prove that $A_\mu^{ss}(\alpha,q)$ is a polynomial in $q$ with rational coefficients.

Let $\mathbb{Q}(q)$ be the field of rational functions in $q$ over $\mathbb{Q}$. The ring of formal power series $R = \mathbb{Q}(q)[[X_1, \dots, X_n]]$ is 
endowed with a $\lambda$-ring structure via the Adams operators $\psi_k$, defined by
\[
\psi_k(f(q, X_1, \dots, X_n)) = f(q^k, X_1^k, \dots, X_n^k).
\]
Let $\mathfrak{m}$ denote the unique maximal ideal of $R$. The \textit{plethystic exponential} map $\operatorname{Exp} : \mathfrak{m} \to 1 + \mathfrak{m}$ is defined by:
\begin{equation} \label{def_Exp}
\operatorname{Exp}(f) = \exp\bigg(\sum_{k \ge 1} \frac{1}{k} \psi_k(f)\bigg).
\end{equation}
This map has an inverse $\operatorname{Log} : 1 + \mathfrak{m} \to \mathfrak{m}$ (see Getzler~\cite{EG 1996} and Mozgovoy~\cite{SM 2007}), called \textit{plethystic logarithm}, given by
\begin{equation} \label{def_Log}
\operatorname{Log}(f) = \sum_{k \ge 1} \frac{\mu(k)}{k} \psi_k(\log(f)),
\end{equation}
where $\mu$ denotes the classical Möbius function.

For $\alpha \in \mathbb{N}^n$, we set $X^\alpha = X_1^{\alpha_1} \cdots X_n^{\alpha_n}$.  The polynomial
$A_\mu^{s}(\alpha, q)$ can be computed using the following identity due to Mozgovoy and Reineke~\cite{M-R 2009}:
\begin{equation} \label{num_of_semistables}
\bigg(\sum_{\alpha\in\Delta_\mu^+} \frac{|\mathrm{R}_\mu^{ss}(\alpha, \mathbb{F}_q)|} {|\mathrm{GL}(\alpha, \mathbb{F}_q)|}  X^\alpha\bigg) \circ 
\text{Exp}\bigg(\frac{1}{1-q}\sum_{\alpha\in\Delta_\mu^+\backslash\{0\}}\!A_\mu^{s}(\alpha, q) X^\alpha\bigg) = 1,
\end{equation}
in the twisted power series ring $\mathbb{Q}(q)^{tw}[[X_1,\ldots,X_n]]$, as defined in Mozgovoy and Reineke~\cite{M-R 2009}.

The main result of this paper is Theorem \ref{main_formula_for_A} which states that 
\begin{equation}\label{main_formula}
A_\mu^{ss}(\alpha,q) = (q-1)\sum_{d \, \mid \, \bar{\alpha}} \frac{\mu(d)}{d}\,H_\mu\!\left(\frac{\alpha}{d}, q^d\right),
\end{equation}
where the sum runs over all divisors of
\[
\bar{\alpha} := \gcd(\alpha_1,\dots,\alpha_n),
\]
and $\mu(d)$ is the classical Möbius function. We refer to Section~3 for the precise definition of
$H_\mu(\alpha,q) \in \mathbb{Q}(q)$.

The generating functions studied in this paper are closely related to the BPS invariants of quivers with stability conditions in the sense of Kontsevich–Soibelman \cite{K-S 2011}. 
In particular, our formula may be interpreted as describing the contribution of semistable representations of fixed slope to the corresponding BPS Lie algebra. 
We do not pursue this perspective here.

This paper is organized as follows:

Section 2 shows that any representation of a quiver over a field admitting a nilpotent endomorphism gives rise to a structured filtration of subrepresentations. 
We then prove that such a representation is semistable of slope $\mu$ if and only if all consecutive quotients in the filtration are semistable of slope $\mu$.
This result allows us to derive the key formula for the number of semistable representations fixed by a given conjugacy class of $\mathrm{GL}(\alpha, \mathbb{F}_q)$.

Section 3 introduces a generating function analogous to that of Hua \cite{JH 2000}, encoding the essential quantities required for the Burnside orbit-counting formula. 
It then establishes a series of identities relating the number of isomorphism classes 
of semistable representations, indecomposable semistable representations, and absolutely indecomposable semistable representations. 
The main result is presented in Theorem~\ref{main_formula_for_A}, with its proof adapted from Hua \cite{JH 2000}, followed by an identity in the style of Mozgovoy~\cite{SM 2007}.

Section 4 introduces refined Kac functions for $A_\mu^{ss}(\alpha, q) \in \mathbb{Q}(q)$ and proposes a conjecture stating that, for quivers of infinite representation type, 
these functions are in fact polynomials in $q$ with non-negative integer coefficients.

Finally, the appendix provides several examples of $A_\mu^{ss}(\alpha, q)$ for quivers with two vertices, offering evidence supporting the conjectures proposed in Section 3 and Section 4.

\section{Filtrations via Nilpotent Endomorphisms}

A \textit{partition} of a non-negative integer $m$ is a sequence of positive integers in non-increasing order whose sum is $m$. The unique partition of $0$ consists of no parts
and is typically represented by the empty sequence $[\,]$. 

Let $\mathcal{P}$ denote the set of all partitions of non-negative integers.
For $\lambda \in \mathcal{P}$, the \textit{weight} of $\lambda$, denoted $|\lambda|$, is the sum of its parts, and 
the \textit{length} of $\lambda$, denoted $l(\lambda)$, is the number of parts in $\lambda$.
For a positive integer $k$,  the \textit{multiplicity} of $k$ in $\lambda$, denoted $m_\lambda^{(k)}$, is defined by
\[
m_\lambda^{(k)} := |\{i : \lambda_i = k \text { where } 1 \leq i \leq l(\lambda) \}|.
\]

For a partition $\lambda=(\lambda_1, \lambda_2, \ldots, \lambda_r) \in \mathcal{P}$, define $J_\lambda$ to be the block-diagonal matrix of size $|\lambda| \times |\lambda|$
given by
\begin{align*}
J_{\lambda} = \left[
\begin{array}{cccc}
J_{\lambda_1} & 0 & \dots & 0 \\ 
0 & \!\!\!\!J_{\lambda_2} & \dots & 0 \\ 
\vdots & \vdots & \ddots & \vdots \\
0 & 0 & \dots & \!\!\!\!J_{\lambda_r}
\end{array}
\right],
\end{align*}
where each $J_{\lambda_i}$ is the $\lambda_i \times \lambda_i$ Jordan block of the form
\[
J_{\lambda_i} = 
\left[
\begin{array}{ccccc}
0 & 1 & 0 & \dots & 0 \\ 
0 & 0 & 1 & \dots & 0 \\ 
\vdots & \vdots & \vdots & \ddots & \vdots \\
0 & 0 & 0 & \dots & 1 \\
0 & 0 & 0 & \dots & 0 
\end{array}
\right].
\]
Every nilpotent matrix $M$ over a field is conjugate to a matrix of the form $J_\lambda$ for a unique partition $\lambda \in \mathcal{P}$. 
The partition $\lambda$ is called the \textit{type} of $M$.

\begin{dfn}
Let $M = [e_{ij}]$ be an $m \times n$ matrix over a field $k$.  
We say that $M$ is a \textit{TA-matrix} if it satisfies the following conditions:
\begin{itemize}
    \item $e_{ij} = e_{kl}$ whenever $j - i = l - k$ for $1\le i,k\le m$ and $1\le j,l\le n$;
    \item $e_{ij} = 0$ if $(i - 1) + (n - j) \geq \min(m, n)$ for $1\le i\le m$ and $1\le j\le n$.
\end{itemize}
\end{dfn}

Thus, a TA-matrix of order $m \times n$ has one of the following forms:
\begin{itemize}
\item
when $m \geq n$, the matrix has the form
\[
\left[
\begin{array}{ccccc}
a_1 & a_2 & \dots & a_{n-1} & a_n \\ 
0 & a_1 & \dots & a_{n-2} & a_{n-1} \\ 
\vdots & \vdots & \ddots & \vdots & \vdots \\
0 & 0 & \dots & a_1 & a_2 \\
0 & 0 & \dots & 0 & a_1 \\
\cline{1-5}
0 & 0 & \dots & 0 & 0 \\
\vdots & \vdots & \vdots & \vdots & \vdots \\
0 & 0 & \dots & 0 & 0
\end{array}
\right],
\]
\item
when $m < n$, it takes the form
\[
\newcommand*{\temp}{\multicolumn{1}{|}{}}
\left[
\begin{array}{ccccccccc}
0 & \dots & 0 & \temp & a_1 & a_2 & \dots & a_{m-1} & a_m \\ 
0 & \dots & 0 & \temp & 0 & a_1 & \dots & a_{m-2} & a_{m-1} \\ 
\vdots & \vdots & \vdots & \temp & \vdots & \vdots & \ddots & \vdots & \vdots \\
0 & \dots & 0 & \temp & 0 & 0 & \dots & a_1 & a_2 \\
0 & \dots & 0 & \temp & 0 & 0 & \dots & 0 & a_1 
\end{array}
\right].
\]
\end{itemize}

\begin{thm}[Turnbull \& Aitken \cite{T-A 1948}]\label{TA Theorem}
Let $\lambda = (\lambda_1, \lambda_2, \dots, \lambda_s)$ and $\mu = (\mu_1, \mu_2, \dots, \mu_t)$ be two partitions.  
A matrix $M$ over a field $k$ satisfies  
\[
J_\lambda M = M J_\mu
\]
if and only if $M$ can be written as an $s \times t$ block matrix of the form:
\[
\left[
\begin{array}{cccc}
T_{11} & T_{12} & \dots & T_{1t} \\ 
T_{21} & T_{22} & \dots & T_{2t} \\ 
\vdots & \vdots & \ddots & \vdots \\
T_{s1} & T_{s2} & \dots & T_{st} 
\end{array}
\right],
\]
where each block $T_{ij}$ is a TA-matrix of size $\lambda_i \times \mu_j$, for $1 \leq i \leq s$ and $1 \leq j \leq t$.
\end{thm}

Let $\lambda \in \mathcal{P}$. Given an integer $r\geq \lambda_1$, let $\eta$ be the partition such that for $1 \le i \le r$,
\[
m_\eta^{(i)} = 
\begin{cases}
	m_\lambda^{(i)} & \text{ if } m_\lambda^{(i)} > 0, \\
	1 &                         \text{ otherwise}.
\end{cases}
\]

The partition $\eta$ defined above is called the \textit{$r$-padded partition} of $\lambda$, denoted by $p_r(\lambda)$. 

The Young diagram for the padded partition will be depicted using 2 colors:
existing parts from $\lambda$ will be shown as black boxes, while padded parts will be shown in blue. 
For example, the Young diagram of $p_4(\lambda)$, where $\lambda=[3,3,1]$, is given below:

\[
\begin{array}{c@{\hskip 3pt}c@{\hskip 3pt}c@{\hskip 3pt}c@{\hskip 3pt}}
\fcolorbox{blue}{white}{\phantom{0}} & \fcolorbox{blue}{white}{\phantom{2}} & \fcolorbox{blue}{white}{\phantom{3}} & \fcolorbox{blue}{white}{\phantom{4}}  \\
[0.1cm]
\boxed{\phantom{1}} & \boxed{\phantom{2}} & \boxed{\phantom{3}} \\
[0.1cm]
\boxed{\phantom{1}} & \boxed{\phantom{2}} & \boxed{\phantom{3}} \\
[0.1cm]
\fcolorbox{blue}{white}{\phantom{1}} & \fcolorbox{blue}{white}{\phantom{2}}  \\
[0.1cm]
\boxed{\phantom{1}} \\
\end{array}
\]

Let $M$ be a representation of the quiver $\Gamma$ over a field $k$, and let $e$ be a nilpotent endomorphism of $M$. A filtration of subrepresentations 
\[
	0 = M_0 \subseteq  M_1 \subseteq \cdots \subseteq M_{r} = M
\]
is said to be $e$-\textit{compatible} if 
\[ 
	e(M_i) \subseteq M_{i-1} \quad \text{for } 1 \leq i \leq r.
\]

\begin{lem}\label{filtration_A2}
Let $k$ be an arbitrary field.  
Given $\lambda, \mu \in \mathcal{P}$ and a matrix $M$ over $k$ such that  
\[
J_\lambda M = M J_\mu,
\]  
we interpret $M$ as a representation of the quiver $A_2$ (i.e., $\bullet \!\! \longrightarrow \!\! \bullet$) over $k$, and the pair $e=(J_\lambda, J_\mu)$ as an endomorphism of this representation.  
Then $M$ admits an $e$-compatible filtration by subrepresentations of the following form:
\[
\begin{array}{rrrrrrrr}
0=M_{0,0} \subseteq &  M_{1,r} \subseteq &  M_{1,r-1} \subseteq & \!\!\!\!\cdots \subseteq & \!\!\!\! M_{1,i} \subseteq & \!\!\!\!\cdots \subseteq & \!\!\!\! M_{1,2} \subseteq & \!\!\!\! M_{1,1}  \\ 
                \subseteq & M_{2,r}  \subseteq & M_{2,r-1}  \subseteq & \!\!\!\!\cdots \subseteq & \!\!\!\! M_{2,i} \subseteq & \!\!\!\!\cdots \subseteq & \!\!\!\! M_{2,2} \phantom{\ \subseteq} \\ 
                              & \vdots \phantom{\subseteq \enspace} & \vdots \phantom{\subseteq \enspace} &                     &\vdots \phantom{\subseteq \enspace} \\
                \subseteq & M_{i,r}   \subseteq & M_{i,r-1}   \subseteq  & \!\!\!\!\cdots \subseteq & \!\!\!\! M_{i,i} \phantom{\ \subseteq}  \\
                              & \vdots \phantom{\subseteq \enspace } & \vdots \phantom{\subseteq \enspace}  \\ 
                \subseteq & \!\!\!\! M_{r-1,r} \subseteq & \!\!\!\! M_{r-1,r-1} \phantom{\ \subseteq} \\
                \subseteq & M_{r,r} = & M, \phantom{\subseteq} \\                
\end{array}
\]
where $r = \max(\lambda_1, \mu_1)$ is the largest part among the partitions $\lambda$ and $\mu$, and the following hold:
\begin{enumerate}
    \item For $1 \le i \le r$ and $i \le s \le r - 1$, we have
    \[
    M_{i,s}/M_{i,s+1} \cong M_{1,s}/M_{1,s+1}, \enspace \text{with} \enspace \dim\left(M_{1,s}/M_{1,s+1}\right) = (m_\lambda^{(s)}, m_\mu^{(s)}).
    \]
    \item For $1 \le i \le r$, we have
    \[
    M_{i,r}/M_{i-1,i-1} \cong M_{1,r}, \enspace \text{with} \enspace \dim M_{1,r} = (m_\lambda^{(r)}, m_\mu^{(r)}).
    \]
\end{enumerate}
\end{lem}
\begin{proof}
Let $V$ and $W$ be vector spaces over $k$, with $\dim V = |\lambda|$ and $\dim W = |\mu|$. Fix a basis $v_1,\ldots, v_m$ for $V$, where $m=\dim V$,
and a basis $w_1, \ldots, w_n$ for $W$, where $n=\dim W$. Define a linear map $f:V \to W$ via the matrix $M$ as follows:
\[
f: \left[
	\begin{array}{c}
		v_1 \\
		v_2 \\ 
		\vdots \\
		v_m \\
	\end{array}
\right] 
\mapsto M
\left[
	\begin{array}{c}
		w_1 \\
		w_2 \\ 
		\vdots \\
		w_n \\
	\end{array}
\right].
\]
Then $(V,W,f)$ defines the representation of  the quiver $A_2$ associated with the matrix $M$. Let $e$ act on $V$ via $J_\lambda$ and on $W$ via $J_\mu$. Thus, 
$e$ is a nilpotent endomorphism of $(V,W,f)$.

Let $s=l(\lambda)$ and $t=l(\mu)$. Since $J_\lambda M = M J_\mu$, Theorem \ref{TA Theorem} implies that $M$ can be expressed as an $s \times t$ block matrix of the form:

\[
\left[
\begin{array}{cccc}
T_{11} & T_{12} & \dots & T_{1t} \\ 
T_{21} & T_{22} & \dots & T_{2t} \\ 
\vdots & \vdots & \ddots & \vdots \\
T_{s1} & T_{s2} & \dots & T_{st} 
\end{array}
\right],
\]
where each block $T_{ij}$ is a TA-matrix of size $\lambda_i \times \mu_j$, for all $1 \leq i \leq s$ and $1 \leq j \leq t$.

We construct a filtration of subrepresentations based on the Young diagrams of $p_r(\lambda)$ and $p_r(\mu)$. 
The basis elements $v_1,v_2,\ldots,v_m$ of $V$ are assigned into the black boxes of the Young diagrams of $p_r(\lambda)$, filled from left to right and top to bottom, while the blue boxes are assigned
the zero vector. Similarly, the basis elements $w_1,w_2,\ldots,w_n$ of $W$ and the zero vector are placed into the Young diagram of $p_r(\mu)$ in the same fashion.

Starting from the zero representation, we progressively select elements from the diagrams to form subspaces of $V$ and $W$, thereby constructing a sequence of subrepresentations of 
$(V,W)$. The procedure is as follows:

\textbf{Step 1.} 
From the Young diagram of $p_r(\lambda)$, select the box at the end of the first row along with all boxes directly beneath it in the same column. The elements associated with these boxes span a subspace of $V$.

Apply the same selection process to the Young diagram of $p_r(\mu)$ to obtain a corresponding subspace of $W$.
Due to the structured pattern of the matrix $M$, the pair of subspaces defines a subrepresentation $M_{1,r}$ of $(V,W)$ with dimension vector $(m_\lambda^{(r)}, m_\mu^{(r)})$. 
The structured pattern of $J_\lambda$ and $J_\mu$ implies that $e(M_{1,r}) = 0$.

Remove the selected boxes from both Young diagrams.

\textbf{Step 2.} 
Identify the first row below the previously selected one in the Young diagram of $p_r(\lambda)$. Select the last box in that row, along with all boxes directly below it in the same column. Form a new subspace of $V$ by taking 
the span of all previously selected elements together with the newly selected ones.

Repeat the same procedure for $W$, using the Young diagram of $p_r(\mu)$. Again the structure of $M$ ensures that the new pair defines a valid subrepresentation $M_{1,r-1}$  of $(V,W)$. Note that the dimension vector of the quotient of this 
representation by the previous one is $(m_\lambda^{(r-1)}, m_\mu^{(r-1)})$. The structured pattern of $J_\lambda$ and $J_\mu$ implies that $e(M_{1,r-1}) \subseteq M_{1,r}$.

Remove the newly selected boxes from both diagrams.
 
\textbf{Step 3.} Continue this process until the last box in the bottom row is selected and $M_{1,1}$ is constructed. The structured pattern of $J_\lambda$ and $J_\mu$ implies that $e(M_{1,1}) \subseteq M_{1,2}$.

Afterward, return to the first row of the remaining diagrams and repeat the procedure until all boxes in the diagrams have been removed. 

At the end of this procedure, the constructed subrepresentation will coincide with the full representation $(V,W)$.

The structured pattern of $M$, $J_\lambda$, and $J_\mu$ ensures that the filtration constructed in this manner is $e$-compatible and satisfies the properties (1) and (2) in the theorem.

Finally, suppose we reorder the basis of $V$ so that the vectors $v_i$ appear in the order in which they were selected, and likewise reorder the basis of $W$
so that the vectors $w_i$ follow their selection order.

Let $U$ be the transition matrix from the original basis of $V$ to the new one, and let $T$ be the transition matrix from the original basis of 
$W$ to its new one. Then, with respect to these new bases, the matrix representing $f$ becomes $UMT^{-1}$. The constructed filtration implies that $UMT^{-1}$ is lower triangular.
Note that both $U$ and $T$ are simply permutation matrices.
\end{proof}

As an example of the above lemma, consider $\lambda=[3,3,2]$ and $\mu=[3,3,1]$. Then we have

\[
J_\lambda =
\left[
\begin{array}{cccccccc}
	0 & 1 & 0 & 0 & 0 & 0 & 0 & 0 \\
	0 & 0 & 1 & 0 & 0 & 0 & 0 & 0 \\
	0 & 0 & 0 & 0 & 0 & 0 & 0 & 0 \\
	0 & 0 & 0 & 0 & 1 & 0 & 0 & 0 \\
	0 & 0 & 0 & 0 & 0 & 1 & 0 & 0 \\
	0 & 0 & 0 & 0 & 0 & 0 & 0 & 0 \\
	0 & 0 & 0 & 0 & 0 & 0 & 0 & 1 \\
	0 & 0 & 0 & 0 & 0 & 0 & 0 & 0 \\
\end{array}
\right], \enspace
J_\mu =
\left[
\begin{array}{cccccccc}
	0 & 1 & 0 & 0 & 0 & 0 & 0 \\
	0 & 0 & 1 & 0 & 0 & 0 & 0 \\
	0 & 0 & 0 & 0 & 0 & 0 & 0 \\
	0 & 0 & 0 & 0 & 1 & 0 & 0 \\
	0 & 0 & 0 & 0 & 0 & 1 & 0 \\
	0 & 0 & 0 & 0 & 0 & 0 & 0 \\
	0 & 0 & 0 & 0 & 0 & 0 & 0 \\
\end{array}
\right].
\]

By Theorem \ref{TA Theorem}, the condition $J_\lambda M = MJ_\mu$ implies that $M$ must have the following structured form:
\[
\left[
\begin{array}{cccccccc}
	a & b & c & d & e & f & g \\
	0 & a & b & 0 & d & e & 0 \\
	0 & 0 & a & 0 & 0 & d & 0 \\
	h & i & j & k & l & m & n \\
	0 & h & i & 0 & k & l & 0 \\
	0 & 0 & h & 0 & 0 & k & 0 \\
	0 & p & q & 0 & r & s & t \\
	0 & 0 & p & 0 & 0 & r & 0 \\
\end{array}
\right].
\]

Furthermore, we have
\[
	p_3(\lambda) = [3,3,2,1] \text{ and } p_3(\mu) = [3,3,2,1], 
\]
and their corresponding Young diagrams are shown below:
\[
\begin{array}{c@{\hskip 3pt}c@{\hskip 3pt}c@{\hskip 3pt}}
\fcolorbox{black}{white}{\phantom{1}} & \fcolorbox{black}{white}{\phantom{1}} & \fcolorbox{black}{white}{\phantom{1}} \\
[0.1cm]
\fcolorbox{black}{white}{\phantom{1}} & \fcolorbox{black}{white}{\phantom{1}} & \fcolorbox{black}{white}{\phantom{1}} \\
[0.1cm]
\fcolorbox{black}{white}{\phantom{1}} & \fcolorbox{black}{white}{\phantom{1}} \\
[0.1cm]
\fcolorbox{blue}{white}{\phantom{1}} \\
\end{array}, \enspace 
\begin{array}{c@{\hskip 3pt}c@{\hskip 3pt}c@{\hskip 3pt}}
\fcolorbox{black}{white}{\phantom{1}} & \fcolorbox{black}{white}{\phantom{1}} & \fcolorbox{black}{white}{\phantom{1}} \\
[0.1cm]
\fcolorbox{black}{white}{\phantom{1}} & \fcolorbox{black}{white}{\phantom{1}} & \fcolorbox{black}{white}{\phantom{1}} \\
[0.1cm]
\fcolorbox{blue}{white}{\phantom{1}} & \fcolorbox{blue}{white}{\phantom{1}} \\
[0.1cm]
\fcolorbox{black}{white}{\phantom{1}} \\
\end{array}.
\]

Let $V$ be an 8-dimensional vector space and $W$ a 7-dimensional vector space, with respective bases
\[
	[v_1, v_2,v_3,v_4,v_5,v_6,v_7,v_8], \, [w_1,w_2,w_3,w_4,w_5,w_6,w_7].
\]

By allocating these basis elements into the Young diagrams of $p_3(\lambda)$ and $p_3(\mu)$, we obtain:

\[
\begin{array}{c@{\hskip 3pt}c@{\hskip 3pt}c@{\hskip 3pt}}
\fcolorbox{black}{white}{$v_1$} & \fcolorbox{black}{white}{$v_2$} & \fcolorbox{black}{white}{$v_3$} \\
[0.1cm]
\fcolorbox{black}{white}{$v_4$} & \fcolorbox{black}{white}{$v_5$} & \fcolorbox{black}{white}{$v_6$} \\
[0.1cm]
\fcolorbox{black}{white}{$v_7$} & \fcolorbox{black}{white}{$v_8$} \\
[0.1cm]
\fcolorbox{blue}{white}{\,0\,\,} \\
\end{array}, \enspace 
\begin{array}{c@{\hskip 3pt}c@{\hskip 3pt}c@{\hskip 3pt}}
\fcolorbox{black}{white}{$w_1$} & \fcolorbox{black}{white}{$w_2$} & \fcolorbox{black}{white}{$w_3$} \\
[0.1cm]
\fcolorbox{black}{white}{$w_4$} & \fcolorbox{black}{white}{$w_5$} & \fcolorbox{black}{white}{$w_6$} \\
[0.1cm]
\fcolorbox{blue}{white}{\,\,0\,\,} & \fcolorbox{blue}{white}{\,\,0\,\,} \\
[0.1cm]
\fcolorbox{black}{white}{$w_7$} \\
\end{array}.
\]

The construction of the filtration of subrepresentations is illustrated in the following diagrams. Rather than removing the selected boxes, we highlight them in red.
\[
\left[
\begin{array}{c@{\hskip 3pt}c@{\hskip 3pt}c@{\hskip 3pt}}
\boxed{v_1} & \boxed{v_2} & \fcolorbox{red}{white}{$v_3$} \\
[0.1cm]
\boxed{v_4} & \boxed{v_5} & \fcolorbox{red}{white}{$v_6$} \\
[0.1cm]
\boxed{v_7} & \boxed{v_8}  \\
[0.1cm]
\fcolorbox{blue}{white}{\,0\,\,} \\
\end{array},
\begin{array}{c@{\hskip 3pt}c@{\hskip 3pt}c@{\hskip 3pt}}
\boxed{w_1} & \boxed{w_2} & \fcolorbox{red}{white}{$w_3$}  \\
[0.1cm]
\boxed{w_4} & \boxed{w_5} & \fcolorbox{red}{white}{$w_6$}  \\
[0.1cm]
\fcolorbox{blue}{white}{\,\,0\,\,}  & \fcolorbox{blue}{white}{\,\,0\,\,} \\
[0.1cm]
\boxed{w_7} \\
\end{array}
\right],
\left[
\begin{array}{c@{\hskip 3pt}c@{\hskip 3pt}c@{\hskip 3pt}}
\boxed{v_1} & \boxed{v_2} & \fcolorbox{red}{white}{$v_3$} \\
[0.1cm]
\boxed{v_4} & \boxed{v_5} & \fcolorbox{red}{white}{$v_6$} \\
[0.1cm]
\boxed{v_7} & \fcolorbox{red}{white}{$v_8$} \\
[0.1cm]
\fcolorbox{blue}{white}{\,0\,\,} \\
\end{array},
\begin{array}{c@{\hskip 3pt}c@{\hskip 3pt}c@{\hskip 3pt}}
\boxed{w_1} & \boxed{w_2} & \fcolorbox{red}{white}{$w_3$}  \\
[0.1cm]
\boxed{w_4} & \boxed{w_5} & \fcolorbox{red}{white}{$w_6$}  \\
[0.1cm]
\fcolorbox{blue}{white}{\,\,0\,\,}  & \fcolorbox{red}{white}{\,\,0\,\,} \\
[0.1cm]
\boxed{w_7} \\
\end{array}
\right],
\]
\[
\left[
\begin{array}{c@{\hskip 3pt}c@{\hskip 3pt}c@{\hskip 3pt}}
\boxed{v_1} & \boxed{v_2} & \fcolorbox{red}{white}{$v_3$} \\
[0.1cm]
\boxed{v_4} & \boxed{v_5} & \fcolorbox{red}{white}{$v_6$} \\
[0.1cm]
\boxed{v_7} & \fcolorbox{red}{white}{$v_8$} \\
[0.1cm]
\fcolorbox{red}{white}{\,0\,\,} \\
\end{array},
\begin{array}{c@{\hskip 3pt}c@{\hskip 3pt}c@{\hskip 3pt}}
\boxed{w_1} & \boxed{w_2} & \fcolorbox{red}{white}{$w_3$}  \\
[0.1cm]
\boxed{w_4} & \boxed{w_5} & \fcolorbox{red}{white}{$w_6$}  \\
[0.1cm]
\fcolorbox{blue}{white}{\,\,0\,\,}  & \fcolorbox{red}{white}{\,\,0\,\,} \\
[0.1cm]
\fcolorbox{red}{white}{$w_7$} \\
\end{array}
\right],
\left[
\begin{array}{c@{\hskip 3pt}c@{\hskip 3pt}c@{\hskip 3pt}}
\boxed{v_1} & \fcolorbox{red}{white}{$v_2$} & \fcolorbox{red}{white}{$v_3$} \\
[0.1cm]
\boxed{v_4} & \fcolorbox{red}{white}{$v_5$} & \fcolorbox{red}{white}{$v_6$} \\
[0.1cm]
\boxed{v_7} & \fcolorbox{red}{white}{$v_8$} \\
[0.1cm]
\fcolorbox{red}{white}{\,0\,\,} \\
\end{array},
\begin{array}{c@{\hskip 3pt}c@{\hskip 3pt}c@{\hskip 3pt}}
\boxed{w_1} & \fcolorbox{red}{white}{$w_2$} & \fcolorbox{red}{white}{$w_3$}  \\
[0.1cm]
\boxed{w_4} & \fcolorbox{red}{white}{$w_5$} & \fcolorbox{red}{white}{$w_6$}  \\
[0.1cm]
\fcolorbox{blue}{white}{\,\,0\,\,}  & \fcolorbox{red}{white}{\,\,0\,\,} \\
[0.1cm]
\fcolorbox{red}{white}{$w_7$} \\
\end{array}
\right],
\]
\[
\left[
\begin{array}{c@{\hskip 3pt}c@{\hskip 3pt}c@{\hskip 3pt}}
\boxed{v_1} & \fcolorbox{red}{white}{$v_2$} & \fcolorbox{red}{white}{$v_3$} \\
[0.1cm]
\boxed{v_4} & \fcolorbox{red}{white}{$v_5$} & \fcolorbox{red}{white}{$v_6$} \\
[0.1cm]
\fcolorbox{red}{white}{$v_7$} & \fcolorbox{red}{white}{$v_8$} \\
[0.1cm]
\fcolorbox{red}{white}{\,0\,\,} \\
\end{array},
\begin{array}{c@{\hskip 3pt}c@{\hskip 3pt}c@{\hskip 3pt}}
\boxed{w_1} & \fcolorbox{red}{white}{$w_2$} & \fcolorbox{red}{white}{$w_3$}  \\
[0.1cm]
\boxed{w_4} & \fcolorbox{red}{white}{$w_5$} & \fcolorbox{red}{white}{$w_6$}  \\
[0.1cm]
\fcolorbox{red}{white}{\,\,0\,\,}  & \fcolorbox{red}{white}{\,\,0\,\,} \\
[0.1cm]
\fcolorbox{red}{white}{$w_7$} \\
\end{array}
\right],
\left[
\begin{array}{c@{\hskip 3pt}c@{\hskip 3pt}c@{\hskip 3pt}}
\fcolorbox{red}{white}{$v_1$} & \fcolorbox{red}{white}{$v_2$} & \fcolorbox{red}{white}{$v_3$} \\
[0.1cm]
\fcolorbox{red}{white}{$v_4$} & \fcolorbox{red}{white}{$v_5$} & \fcolorbox{red}{white}{$v_6$} \\
[0.1cm]
\fcolorbox{red}{white}{$v_7$} & \fcolorbox{red}{white}{$v_8$} \\
[0.1cm]
\fcolorbox{red}{white}{\,0\,\,} \\
\end{array},
\begin{array}{c@{\hskip 3pt}c@{\hskip 3pt}c@{\hskip 3pt}}
\fcolorbox{red}{white}{$w_1$} & \fcolorbox{red}{white}{$w_2$} & \fcolorbox{red}{white}{$w_3$}  \\
[0.1cm]
\fcolorbox{red}{white}{$w_4$} & \fcolorbox{red}{white}{$w_5$} & \fcolorbox{red}{white}{$w_6$}  \\
[0.1cm]
\fcolorbox{red}{white}{\,\,0\,\,}  & \fcolorbox{red}{white}{\,\,0\,\,} \\
[0.1cm]
\fcolorbox{red}{white}{$w_7$} \\
\end{array}
\right],
\]

Thus, we obtain the following filtration of subrepresentations:
\begin{align*}
(0,0) & \subseteq (\langle v_3, v_6\rangle, \langle w_3, w_6\rangle) \subseteq (\langle v_3, v_6, v_8\rangle, \langle w_3, w_6\rangle) \subseteq (\langle v_3, v_6, v_8\rangle, \langle w_3, w_6, w_7\rangle) \\
& \subseteq (\langle v_3, v_6, v_8,v_2,v_5\rangle, \langle w_3, w_6, w_7,w_2,w_5\rangle)  \\
& \subseteq (\langle v_3, v_6, v_8,v_2,v_5,v_7\rangle, \langle w_3, w_6, w_7,w_2,w_5\rangle)  \\
& \subseteq (\langle v_3, v_6, v_8,v_2,v_5,v_7,v_1,v_4\rangle, \langle w_3, w_6, w_7,w_2,w_5,w_1,w_4\rangle) = (V,W)
\end{align*}
Let $U$ be the transition matrix that maps
\[
[v_1,v_2,v_3,v_4,v_5,v_6,v_7,v_8] \mapsto [v_3, v_6, v_8,v_2,v_5,v_7,v_1,v_4],
\]
and let $T$ be the transition matrix that maps
\[
[w_1,w_2,w_3,w_4,w_5,w_6,w_7] \mapsto [w_3,w_6,w_7,w_2,w_5,w_1,w_4].
\]
Then we have
\[
U =
\left[
\begin{array}{cccccccc}
	0 & 0 & 1 & 0 & 0 & 0 & 0 & 0 \\
	0 & 0 & 0 & 0 & 0 & 1 & 0 & 0 \\
	0 & 0 & 0 & 0 & 0 & 0 & 0 & 1 \\
	0 & 1 & 0 & 0 & 0 & 0 & 0 & 0 \\
	0 & 0 & 0 & 0 & 1 & 0 & 0 & 0 \\
	0 & 0 & 0 & 0 & 0 & 0 & 1 & 0 \\
	1 & 0 & 0 & 0 & 0 & 0 & 0 & 0 \\
	0 & 0 & 0 & 1 & 0 & 0 & 0 & 0 \\
\end{array}
\right], \enspace
T =
\left[
\begin{array}{cccccccc}
	0 & 0 & 1 & 0 & 0 & 0 & 0 \\
	0 & 0 & 0 & 0 & 0 & 1 & 0 \\
	0 & 0 & 0 & 0 & 0 & 0 & 1 \\
	0 & 1 & 0 & 0 & 0 & 0 & 0 \\
	0 & 0 & 0 & 0 & 1 & 0 & 0 \\
	1 & 0 & 0 & 0 & 0 & 0 & 0 \\
	0 & 0 & 0 & 1 & 0 & 0 & 0 \\
\end{array}
\right].
\]
Simple calculation shows that
\[
UMT^{-1} = 
\left[
\begin{array}{cccccccc}
	a & d & 0 & 0 & 0 & 0 & 0 \\
	h & k & 0 & 0 & 0 & 0 & 0 \\
	p & r & 0 & 0 & 0 & 0 & 0 \\
	b & e & 0 & a & d & 0 & 0 \\
	i & l & 0 & h & k & 0 & 0 \\
	q & s & t & p & r & 0 & 0 \\
	c & f & g & b & e & a & d \\
	j & m & n & i & l & h & k \\
\end{array}
\right],
\]
which is a lower triangular matrix, as expected. The pattern of the matrix above confirms the properties claimed by the lemma.

For an $n$-tuple of partitions $\lambda_* = (\lambda^1, \dots, \lambda^n) \in \mathcal{P}^n$ and an integer $k \geq 1$, the \textit{multiplicity vector} of $k$ in $\lambda_*$, denoted $\lambda_*^{(k)}$,
is defined by
\begin{equation} \label{k_multiplicity_vector}
	\lambda_*^{(k)} = \big( m_{\lambda^1}^{(k)}, \dots, m_{\lambda^n}^{(k)} \big).
\end{equation}

Lemma~\ref{filtration_A2}, together with its proof, extends naturally to an arbitrary quiver.

\begin{thm}\label{filtration_for_quiver}
Let $M$ be a representation of the quiver $\Gamma$ over a field $k$, and let $e=(e_1,\ldots,e_n)$ be a nilpotent endomorphism of $M$. 
Assume that $e_i$ has type $\lambda^i \in \mathcal{P}$ for $1\le i \le n$ and let $\lambda_*= (\lambda^1,\dots,\lambda^n) \in \mathcal{P}^n$. 
Then $M$ admits an $e$-compatible filtration by subrepresentations of the following form:
\[
\begin{array}{rrrrrrrr}
0=M_{0,0} \subseteq &  M_{1,r} \subseteq &  M_{1,r-1} \subseteq & \!\!\!\!\cdots \subseteq & \!\!\!\! M_{1,i} \subseteq & \!\!\!\!\cdots \subseteq & \!\!\!\! M_{1,2} \subseteq & \!\!\!\! M_{1,1}  \\ 
                \subseteq & M_{2,r}  \subseteq & M_{2,r-1}  \subseteq & \!\!\!\!\cdots \subseteq & \!\!\!\! M_{2,i} \subseteq & \!\!\!\!\cdots \subseteq & \!\!\!\! M_{2,2} \phantom{\ \subseteq} \\ 
                              & \vdots \phantom{\subseteq \enspace} & \vdots \phantom{\subseteq \enspace} &                     &\vdots \phantom{\subseteq \enspace} \\
                \subseteq & M_{i,r}   \subseteq & M_{i,r-1}   \subseteq  & \!\!\!\!\cdots \subseteq & \!\!\!\! M_{i,i} \phantom{\ \subseteq}  \\
                              & \vdots \phantom{\subseteq \enspace } & \vdots \phantom{\subseteq \enspace}  \\ 
                \subseteq & \!\!\!\! M_{r-1,r} \subseteq & \!\!\!\! M_{r-1,r-1} \phantom{\ \subseteq} \\
                \subseteq & M_{r,r} = & M, \phantom{\subseteq} \\                
\end{array}
\]
where $r$ is the largest part among the partitions $\lambda^1, \ldots, \lambda^n$, and the following hold:
\begin{enumerate}
    \item For $1 \le i \le r$ and $i \le s \le r - 1$, we have
    \[
    M_{i,s}/M_{i,s+1} \cong M_{1,s}/M_{1,s+1}, \enspace \text{with} \enspace \dim\left(M_{1,s}/M_{1,s+1}\right) = \lambda_*^{(s)}.
    \]
    \item For $1 \le i \le r$, we have
    \[
    M_{i,r}/M_{i-1,i-1} \cong M_{1,r}, \enspace \text{with} \enspace \dim M_{1,r} = \lambda_*^{(r)}.
    \]
\end{enumerate}
\end{thm}

A key observation is that Theorem~\ref{filtration_for_quiver} applies to the $g$-loop quiver for $g\geq 1$, which has a single vertex and $g$ edge-loops. Thus,  
it also applies to any finitely generated algebra over an arbitrary field.

The next two lemmas appear in Reineke \cite{M-R 2008}, and are included here for convenience.
\begin{lem} \label{lem_R1}
Let $0 \to X \to Y \to Z \to 0$ be a short exact sequence of non-zero representations of $\Gamma$ over a field $k$. 
Then we have
\begin{enumerate}
  \item $\mu(X) \le \mu(Y)$ if and only if $\mu(Y) \le \mu(Z)$ if and only if $\mu(X) \le \mu(Z)$.
  \item $\mu(X) < \mu(Y)$ if and only if $\mu(Y) < \mu(Z)$ if and only if $\mu(X) < \mu(Z)$.
\end{enumerate}
\end{lem}

\begin{lem} \label{lem_R2}
Let $0 \to X \to Y \to Z \to 0$ be a short exact sequence of non-zero representations of $\Gamma$ over a field $k$ and
$\mu(X) = \mu(Y) = \mu(Z)$. Then $Y$ is semistable if and only if $X$ and $Z$ are semistable.
\end{lem}

\begin{lem} \label{lem_H1}
Let $0 \to X \to Y \to Z \to 0$ be a short exact sequence of non-zero representations of $\Gamma$ over a field $k$. 
If $\mu(X) \le \mu(Y)$ and $\mu(Z) \le \mu(Y)$, then $\mu(X) = \mu(Y) = \mu(Z)$.
\end{lem}

\begin{proof}
$\mu(X) \le \mu(Y)$ implies that $\mu(Y) \le \mu(Z)$ by Lemma \ref{lem_R1}, thus $\mu(Z) = \mu(Y)$. Similarly, $\mu(X) = \mu(Y)$.
\end{proof}

\begin{cor} \label{semistable iff}
Let $M$ and $\lambda$ be as in Theorem~\ref{filtration_for_quiver}, and assume that $M_{1,r+1} = 0$. Then $M$ is semistable of slope $\mu$ if and only if, for all $1 \le s \le r$, 
the quotient $M_{1,s}/M_{1,s+1}$ is either zero or semistable of slope $\mu$.
\end{cor}
\begin{proof}
If $M_{1,s}/M_{1,s+1}$ is either zero or semistable of slope $\mu$ for all $1 \leq s \leq r$, then $M$ is semistable of slope $\mu$, as the category 
of semistable representations of slope $\mu$ is closed under extensions by Lemma~\ref{lem_R2}.

Now suppose that $M$ is semistable and $\mu(M) = \mu$. We have the following short exact sequence:
\[
	0 \to M_{r-1,r-1} \to M \to M_{1,r} \to 0.
\]

If $M_{r-1,r-1} = 0$, then $M \cong M_{1,r}$ and hence $M_{1,r}$ is semistable of slope $\mu$. 
Assume now that $M_{r-1,r-1} \ne 0$. Since $M_{r-1,r-1}$ and $M_{1,r}$ are proper subrepresentations of $M$, we have $\mu(M_{r-1,r-1}) \leq \mu(M)$ 
and $\mu(M_{1,r}) \leq \mu(M)$. 

By Lemma~\ref{lem_H1}, it follows that $\mu(M_{r-1,r-1}) = \mu(M_{1,r}) = \mu$.
Applying Lemma~\ref{lem_R2}, we conclude that both $M_{r-1,r-1}$ and $M_{1,r}$ are semistable.

Suppose that $M_{1,r-1}/M_{1,r} \neq 0$. We have the following short exact sequence:
\[
0 \to M_{1,r} \to M_{1,r-1} \to M_{1,r-1}/M_{1,r} \to 0.
\]
Since $\mu(M_{1,r}) = \mu$ and $\mu(M_{1,r-1}) \leq \mu$, Lemma~\ref{lem_R1} implies that $\mu(M_{1,r-1}/M_{1,r}) \leq \mu$.  
Next, consider the short exact sequence:
\[
0 \to M_{r-1,r} \to M_{r-1,r-1} \to M_{1,r-1}/M_{1,r} \to 0.
\]
By Lemma~\ref{lem_H1}, we have $\mu(M_{r-1,r}) = \mu(M_{1,r-1}/M_{1,r}) = \mu$.  
Applying Lemma~\ref{lem_R2}, it follows that both $M_{r-1,r}$ and $M_{1,r-1}/M_{1,r}$ are semistable.  

An induction argument completes the proof.
\end{proof}

For $\lambda \in \mathcal{P}$, let $\lambda'=(\lambda'_1, \lambda'_2, \ldots)$ denote the \textit{conjugate partition} of $\lambda$, defined by 
\[
	\lambda'_i = \sum_{k\ge i}m_\lambda^{(k)}. 
\]
Given two partitions $\lambda,\mu\in\mathcal{P}$, we define the following pairings:
\begin{align*}
	\langle \lambda, \mu \rangle = \sum_{i\ge 1} \lambda'_i\mu'_i, \quad
	\langle\!\langle \lambda, \mu \rangle\!\rangle = \sum_{i\ge 1} \lambda'_i\mu'_i - \sum_{i\ge 1} m_\lambda^{(i)}m_\mu^{(i)}.
\end{align*}

\begin{lem}\label{alt def}
For $\lambda, \mu \in \mathcal{P}$, we have 
\[
	\langle \lambda, \mu \rangle = \sum_{i=1}^{l(\lambda)} \sum_{j=1}^{l(\mu)} \min(\lambda_i, \mu_j).
\]
\end{lem}
\begin{proof}
Using the notation introduced earlier, we compute:
\begin{align*}
	\sum_{i=1}^{l(\lambda)} \sum_{j=1}^{l(\mu)} \min(\lambda_i, \mu_j) 
	&= \sum_{i=1}^{\lambda_1} \sum_{j=1}^{\mu_1} m_\lambda^{(i)} m_\mu^{(j)} \min(i, j).
\end{align*}
This identity is established in Lemma 3.3 of Hua~\cite{JH 2000}, which shows that
\[
	\sum_{i=1}^{\lambda_1} \sum_{j=1}^{\mu_1} m_\lambda^{(i)} m_\mu^{(j)} \min(i, j) = \langle \lambda, \mu \rangle.
\]
The result follows immediately.
\end{proof}

Conceptually, the difference $\langle \lambda, \mu \rangle - \langle\!\langle \lambda, \mu \rangle\!\rangle$ equals the total number of entries contained in the first $r$ diagonal blocks 
of the matrix $UMT^{-1}$, corresponding to the direct sum $M_{1,r} / M_{1,r+1} \oplus \cdots \oplus M_{1,1} / M_{1,2}$,  as constructed in the proof of Lemma \ref{filtration_A2}.

The following is a direct consequence of Lemma~\ref{filtration_A2} and Lemma~\ref{alt def}.

\begin{lem} \label{num commutators}
Let \( \lambda, \mu \in \mathcal{P} \), and let \( s = |\lambda| \) and \( t = |\mu| \). Then we have
\[
\left| \left\{ M \in \mathrm{Mat}_{s \times t}(\mathbb{F}_q) \,\middle|\, J_\lambda M = M J_\mu \right\} \right| = q^{\langle \lambda, \mu \rangle}.
\]
\end{lem}

\begin{thm}\label{num_of_fixed_rep}
Let $\lambda_* = (\lambda^1, \ldots, \lambda^n) \in \mathcal{P}^n$ be an $n$-tupe of partitions, and define $\alpha=(|\lambda^1|,\ldots,|\lambda^n|) \in \mathbb{N}^n$. Let
\[
	g=(I+J_{\lambda^1},\ldots,I+J_{\lambda^n})\in\textup{GL}(\alpha, \mathbb{F}_q),
\]
where $I$ denotes the identity matrix. Consider the fixed-point set
\[
X_g=\left\{M\in\mathrm{R}_\mu^{ss}(\alpha,\mathbb{F}_q) \mid g\cdot M = M\right\}.
\]
Then, we have
\begin{equation} \label{num_of_fixed_points}
	|X_g| = q^{\sum_{1\le i,j\le n} a_{ij} \langle\!\langle \lambda^i,\lambda^j \rangle\!\rangle} 
	\prod_{s\ge 1} |\mathrm{R}_\mu^{ss}(\lambda_*^{(s)} , \mathbb{F}_q)|.
\end{equation}
\end{thm}
\begin{proof}
Let $M\in\mathrm{Rep}(\alpha,\mathbb{F}_q)$. Then $g\cdot M = M$ if and only if $(J_{\lambda^1}, \ldots, J_{\lambda^n})$ defines an endomorphism of $M$.
By Theorem~\ref{TA Theorem} and Lemma~\ref{num commutators}, it follows that
\[
	|\left\{M\in\mathrm{Rep}(\alpha,\mathbb{F}_q) \mid g\cdot M = M\right\}| =  q^{\sum_{1\le i,j\le n} a_{ij}\langle \lambda^i,\lambda^j\rangle}.
\]

Since each $J_{\lambda^i}$ is nilpotent, Theorem~\ref{filtration_for_quiver} implies that $M$ admits a filtration of subrepresentations with prescribed properties.
By Corollary~\ref{semistable iff}, $M$ is semistable of slope $\mu$ if and only if each quotient $M_{1,s}/M_{1,s+1}$ is either zero or semistable of slope $\mu$ for all $1\le s \le r$.
Moreover, the dimension vector of each quotient $M_{1,s} / M_{1,s+1}$ is $\lambda_*^{(s)}$.

Now consider the representation $M$ in terms of matrices. By Lemma~\ref{filtration_A2} and Theorem~\ref{filtration_for_quiver}, there exist a tuple of permutation matrices over $\mathbb{F}_q$,  
$(U_i)_{i\in\Gamma_0}$, such that for each arrow $i \to j$ in $\Gamma_1$, and for the linear map associated with this arrow in $M$, denoted $K$,  
the matrix $U_i K U_j^{-1}$ is lower triangular. 

Moreover, the first $r$ diagonal blocks correspond to the representation
\[
M_{1,r} / M_{1,r+1} \oplus \cdots \oplus M_{1,1} / M_{1,2}.
\]

The total number of entries contained in the first $r$ blocks is given by
\[
	\sum_{1\le i,j\le n} a_{ij}(\langle \lambda^i,\lambda^j \rangle - \langle\!\langle \lambda^i,\lambda^j \rangle\!\rangle).
\]
Therefore, the number of possible choices for the first $r$ blocks such that
\begin{equation} \label{first r blocks}
M_{1,r} / M_{1,r+1} \oplus \cdots \oplus M_{1,1} / M_{1,2}
\end{equation}
forms a representation is
\[
q^{\sum_{1 \leq i,j \leq n} a_{ij} \left( \langle \lambda^i, \lambda^j \rangle - \langle\!\langle \lambda^i, \lambda^j \rangle\!\rangle \right)}.
\]
Furthermore, the number of choices such that \eqref{first r blocks} forms a semistable representation of slope $\mu$ is
\[
\prod_{s \geq 1} | \mathrm{R}_\mu^{\mathrm{ss}}(\lambda_*^{(s)}, \mathbb{F}_q) |.
\]
Thus, we have
\begin{align*}
	|X_g| 
	&= q^{\sum_{1\le i,j\le n} a_{ij}(\langle \lambda^i,\lambda^j \rangle - (\langle \lambda^i,\lambda^j \rangle - \langle\!\langle \lambda^i,\lambda^j \rangle\!\rangle))} \prod_{s\ge 1} |\mathrm{R}_\mu^{ss}(\lambda_*^{(s)} , \mathbb{F}_q)| \\
	&= q^{\sum_{1\le i,j\le n} a_{ij}\langle\!\langle \lambda^i,\lambda^j \rangle\!\rangle} \prod_{s\ge 1} |\mathrm{R}_\mu^{ss}(\lambda_*^{(s)} , \mathbb{F}_q)| ,
\end{align*}
\end{proof}

\section{Generating Functions for Numbers of Isomorphism Classes}

For $\alpha \in \Delta_\mu^+$, let $M_\mu^{\mathrm{ss}}(\alpha, q)$ (respectively, $I_\mu^{\mathrm{ss}}(\alpha, q)$) be the number of isomorphism classes of semistable 
(respectively, indecomposable semistable) representations of $\Gamma$ with dimension $\alpha$.

The finite group $\mathrm{GL}(\alpha, \mathbb{F}_q)$ acts on the finite set $\mathrm{R}_\mu^{\mathrm{ss}}(\alpha, \mathbb{F}_q)$, and the quantity 
$M_\mu^{\mathrm{ss}}(\alpha, q)$ equals the number of orbits in $\mathrm{R}_\mu^{\mathrm{ss}}(\alpha, \mathbb{F}_q)$.  
Hence, by the Burnside's orbit-counting formula, we have
\begin{equation} \label{oribit_counting}
    M_\mu^{\mathrm{ss}}(\alpha, q) 
    = \frac{1}{|\mathrm{GL}(\alpha, \mathbb{F}_q)|} \sum_{g \in \mathrm{GL}(\alpha, \mathbb{F}_q)} |X_g| 
    = \sum_{g \in \mathrm{CL}(\alpha, \mathbb{F}_q)} \frac{|X_g|}{|\,Z_g|},
\end{equation}
where $\mathrm{CL}(\alpha, \mathbb{F}_q)$ denotes a complete set of representatives of the conjugacy classes of $\mathrm{GL}(\alpha, \mathbb{F}_q)$, and:
\[
    X_g = \{ M \in \mathrm{R}_\mu^{\mathrm{ss}}(\alpha, \mathbb{F}_q) \mid g \cdot M = M \},
\]
\[
    Z_g = \{ h \in \mathrm{GL}(\alpha, \mathbb{F}_q) \mid g^{-1} h g = h \}.
\]

For $\lambda \in \mathcal{P}$, let $b_\lambda(q)$ denote the product as follows:
\[
	b_\lambda(q) = \prod_{i \ge 1} \prod_{k = 1}^{m_\lambda^{(i)}} (1 - q^k).
\]

Let \( g \) be as defined in Theorem~\ref{num_of_fixed_rep}. Then \( |X_g| \) is given by~\eqref{num_of_fixed_points}, 
and \( |Z_g| \) is given by Hua~\cite{JH 2000} as follows:
\[
|Z_g| = \prod_{i=1}^n q^{\langle \lambda^i, \lambda^i \rangle} b_{\lambda^i}(q^{-1}).
\]

\begin{dfn}\label{def_P}
For $\lambda_* = (\lambda^1, \ldots, \lambda^n) \in \mathcal{P}^n$, let $X^{|\lambda_*|} = X_1^{|\lambda^1|} \cdots X_n^{|\lambda^n|}$.
We define a formal power series $P_\mu$ in $\mathbb{Q}(q)[[X_1,\ldots,X_n]]$ by
\begin{equation}\label{def_P_mu}
P_\mu(q, X_1,\ldots,X_n) := \sum_{\lambda_* \in \mathcal{P}^n} \!
\frac{q^{\sum_{1\le i, j \le n}\!a_{ij}\langle\!\langle \lambda^i,\lambda^j \rangle\!\rangle}\!\prod_{s\ge 1}\!|\mathrm{R}_\mu^{ss}(\lambda_*^{(s)}\!, \mathbb{F}_q)|}
{\prod_{1\le i \le n}\!q^{\langle \lambda^i, \lambda^i\rangle} b_{\lambda^i}(q^{-1})} X^{|\lambda_*|},
\end{equation}
Note that the $n$-tuple $([\,], \ldots, [\,]) \in \mathcal{P}^n$ contributes a term equal to $1$ in the sum.
\end{dfn}

The following results are in parallel with those of Hua~\cite{JH 2000}.
\begin{thm} \label{MIA}
With notation as above, we have
\begin{align} 
\sum_{\alpha \in \Delta_\mu^+} M_\mu^{ss}(\alpha, q) X^\alpha &= \prod_{d = 1}^\infty \left(P_\mu(q^d,X_1^d, \dots, X_n^d)\right)^{\phi_d(q)}, \label{burn}\\
\sum_{\alpha \in \Delta_\mu^+} M_\mu^{ss}(\alpha, q) X^\alpha &= \prod_{\alpha \in \Delta_\mu^+ \setminus \{0\}} (1 - X^\alpha)^{-I_\mu^{ss}(\alpha, q)}, \label{KS}\\
I_\mu^{ss}(\alpha, q) &= \sum_{d \,|\, \bar{\alpha}} \frac{1}{d} \sum_{r \,|\, d} \mu \Big( \frac{d}{r} \Big)\, A_\mu^{ss}\!\left(\frac{\alpha}{d}, q^r\right), \label{I_in_terms_of_A}\\
A_\mu^{ss}(\alpha, q) &= \sum_{d \,|\, \bar{\alpha}} \frac{1}{d} \sum_{r \,|\, d} \mu(r)\, I_\mu^{ss}\!\left(\frac{\alpha}{d}, q^r\right), \label{A_in_terms_of_I}
\end{align}
where $\phi_d(q)$ denotes the number of monic irreducible polynomials of degree $d$ in $\mathbb{F}_q[t]$, excluding $t$, and
$\mu(r)$ is the classical Möbius function. Here, $\bar{\alpha} = \gcd(\alpha_1, \dots, \alpha_n)$ denotes the greatest common divisor of $\alpha_1, \dots, \alpha_n$.
\end{thm}

\begin{proof}
Identity~\eqref{burn} follows from the Burnside orbit-counting formula \eqref{oribit_counting} and Theorem~\ref{num_of_fixed_rep}. The same argument as in Theorem 4.3 of Hua~\cite{JH 2000} applies here.

The category of all semistable representations of $\Gamma$ of slope~$\mu$ is a Krull--Schmidt category. 
Hence, every semistable representation of $\Gamma$ of slope~$\mu$ can be uniquely decomposed (up to the order of summands) as a direct sum of indecomposable 
semistable representations of the same slope. Identity~\eqref{KS} follows immediately.

The origin of identity~\eqref{I_in_terms_of_A} can be traced back to Kac~\cite{VK 1983}, and remains valid in the present context.

Identity~\eqref{A_in_terms_of_I} is analogous to the second identity in Theorem~4.1 of Hua~\cite{JH 2000}, and serves as the Möbius inverse of identity~\eqref{I_in_terms_of_A}.
\end{proof}

\begin{dfn}\label{def of H}
We define rational functions $H_\mu(\alpha,q) \in \mathbb{Q}(q)$ for all $\alpha \in \mathbb{N}^n\backslash\{0\}$ as follows:
\begin{equation}\label{log(P)=sum(H)}
\log\left(P_\mu(q, X_1,\ldots,X_n) \right) = \sum_{\alpha\in\mathbb{N}^n\backslash\{0\}} \!H_\mu(\alpha,q)X^\alpha ,
\end{equation}
where $\log$ is the formal logarithm, i.e., $\log(1+x) = \sum_{i\ge 1} (-1)^{i-1} x^i/i$.
\end{dfn}

It follows from \eqref{def_P_mu} that $H_\mu(\alpha, q) = 0$ if $\alpha \notin \Delta_\mu^+$.
The main result of this paper is the following theorem, whose proof is analogous to that of Theorem~4.6 in Hua~\cite{JH 2000}.

\begin{thm}\label{main_formula_for_A} For any $\alpha \in \Delta_\mu^+\setminus \{0\}$, we have:
\begin{equation}\label{A=sum(H)}
A_\mu^{ss}(\alpha,q) = (q-1)\sum_{d \,\mid\, \bar{\alpha}}\frac{\mu(d)}{d}H_\mu\Big(\frac{\alpha}{d}, q^d\Big),
\end{equation}
where the sum runs over all divisors of $\bar{\alpha}=\gcd(\alpha_1, \dots, \alpha_n)$.
\end{thm}

Theorems \ref{MIA} and \ref{main_formula_for_A} imply that $M_\mu^{ss}(\alpha,q)$, $I_\mu^{ss}(\alpha,q)$, and $A_\mu^{ss}(\alpha,q)$ are rational functions in $q$.
Since each of these functions takes integer values for any prime power $q$, they must, in fact, be polynomials in $q$, i.e., 
\[
	M_\mu^{ss}(\alpha, q), I_\mu^{ss}(\alpha, q), A_\mu^{ss}(\alpha, q) \in \mathbb{Q}[q].
\]

Each $n$-tuple of partitions $\lambda_* = (\lambda^1, \ldots, \lambda^n) \in \mathcal{P}^n$ defines an $r$-tuple of dimension vectors 
$\alpha_* = (\alpha^1, \ldots, \alpha^r)$ with $\alpha^k \in \mathbb{N}^n$, where $r$ is the largest part among the partitions 
$\lambda^1, \ldots, \lambda^n$. Here, $\alpha^k$ is the $k$-multiplicity vector of $\lambda_*$, as defined in \eqref{k_multiplicity_vector}, 
and the correspondence is bijective. Replacing $\lambda_*$ with $\alpha_*$ in \eqref{def_P_mu} yields the following equivalent expression for $P_\mu(q, X_1, \ldots, X_n)$:
\begin{equation}\label{def PX}
P_\mu(q, X_1,\ldots,X_n) = 1 + \sum_{\alpha_*} \!
\Big(\prod_{s \ge 1} \frac{q^{\langle \alpha^s, \alpha^s\rangle}} {q^{\langle \beta^s, \beta^s \rangle}}  \frac{|R_\mu^{ss}(\alpha^s, \mathbb{F}_q)|}{|\mathrm{GL}(\alpha^s, \mathbb{F}_q)|} X^{s\alpha^s}\Big),
\end{equation}
where the sum runs over all $r$-tuples of dimension vectors  $\alpha_*=(\alpha^1,\alpha^2,\ldots, \alpha^r)$ with $r \ge 1$, $\alpha^i \in \Delta_\mu^+$ for $1\le i \le r$, and $\alpha^r \ne 0$.
Here $\beta^k=\sum_{i \geq k} \alpha^i$ for $1 \le k \le r$, and $\langle-,-\rangle$ denotes the Euler form associated with the quiver $\Gamma$, as defined in \eqref{euler_form}.

\begin{cor}
The following identity holds in the formal power series ring $\mathbb{Q}(q)[[X_1,\ldots,X_n]]$:
\begin{equation}\label{P=Exp(A/(q-1))}
	P_\mu(q,X_1,\ldots,X_n) = \mathrm{Exp}\Big(\frac{1}{q-1} \sum_{\alpha\in\Delta_\mu^+\backslash\{0\}}\!A_\mu^{ss}(\alpha,q)X^\alpha\Big),
\end{equation}
where $\mathrm{Exp}$ denotes the plethystic exponential map defined in \eqref{def_Exp}.
\end{cor}
\begin{proof}
By the definition of the plethystic logarithm, we have9i9
\begin{align*}
\mathrm{Log} (P_\mu(q,X_1,\ldots,X_n)) &= \sum_{k=1}^\infty \frac{\mu(k)}{k} \log ( P_\mu(q^k,X_1^k,\ldots,X_n^k) ) \\
&\overset{\eqref{log(P)=sum(H)}}{=}\sum_{k=1}^\infty \frac{\mu(k)}{k} \sum_{\alpha\in\mathbb{N}^n\backslash\{0\}} \!H_\mu(\alpha,q^k)X^{k\alpha} \\
&= \sum_{\alpha\in\mathbb{N}^n\backslash\{0\}} \bigg(\sum_{d\,|\,\bar{\alpha}}\frac{\mu(d)}{d} H_\mu\Big(\frac{\alpha}{d}, q^d\Big) \bigg)X^\alpha \\
&\overset{\eqref{A=sum(H)}}{=} \frac{1}{q-1}  \sum_{\alpha\in\Delta_\mu^+ \setminus \{0\}} A_\mu^{ss}(\alpha,q)X^\alpha.
\end{align*}
Taking the plethystic exponential of both sides gives the desired identity:
\[
P_\mu(q,X_1,\ldots,X_n) = \mathrm{Exp}\Big(\frac{1}{q-1} \sum_{\alpha\in\Delta_\mu^+ \setminus \{0\}}\!A_\mu^{ss}(\alpha,q)X^\alpha\Big).
\]
\end{proof}

We conjecture that the polynomial $A_\mu^{\mathrm{ss}}(\alpha, q)$ exhibits a positivity property, meaning that all coefficients of $A_\mu^{\mathrm{ss}}(\alpha, q)$ 
are non-negative integers for each $\alpha\in\Delta_\mu^+ \setminus \{0\}$; see the appendix for examples. It is conceivable that the constant term of $A_\mu^{\mathrm{ss}}(\alpha, q)$ 
encodes meaningful geometric information, possibly analogous to that predicted by the Kac conjecture in \cite{VK 1983}.

\section{Refined Kac Functions for Semistable Representations}

The refined Kac functions were first introduced by Rodriguez-Villegas~\cite{FRV 2011} for the $g$-loop quiver and were later extended to general quivers by Hua~\cite{JH 2026}. 
Hua~\cite{JH 2026} also showed that if the quiver contains at least one edge-loop for each vertex, then the refined Kac functions are polynomials in $q$ with non-negative integer coefficients.

For an $n$-tuple of partitions $\lambda_* = (\lambda^1, \ldots, \lambda^n) \in \mathcal{P}^n$, the vector 
\[
	\lambda_*^{(k)} = \big(m_{\lambda^1}^{(k)}, \ldots, m_{\lambda^n}^{(k)}\big) \in \mathbb{N}^n
\]
is called the multiplicity vector of $k$ in $\lambda_*$. We define an $r$-tupe of vectors as follows:
\[
	\lambda_*^\mathrm{v} := \big(\lambda_*^{(1)}, \ldots, \lambda_*^{(r)}\big),
\]
where $r$ is the maximal part in $\lambda^i$ for $1 \le i \le n$. 

Let $\mathbb{Q}(q)[[X_{ik}]]_{1 \le i \le n, k \ge 1}$ denote the formal power series ring over $\mathbb{Q}(q)$ in the infinitely many variables $X_{ik}$, indexed by $1 \le i \le n$ and $k \ge 1$.
For an $r$-tuple of dimension vectors $\alpha_* = (\alpha^1, \ldots, \alpha^r)$ where $\alpha^k = (\alpha^k_1, \ldots, \alpha^k_n) \in \mathbb{N}^n$ for $1 \le k \le r$, we define
\[
	X^{\alpha_*} = \prod_{i=1}^n \prod_{k=1}^r X_{ik}^{k \alpha^k_i}.
\]

Following Hua \cite{JH 2026}, we define a formal power series $Q_\mu$ in $\mathbb{Q}(q)[[X_{ik}]]$ as follows:
\begin{equation}\label{def_Q}
Q_\mu(q, X_{ik})_{1\le i \le n, k\ge 1}:= 1 + \sum_{\alpha_*} \!
\bigg(\prod_{k=1}^r \frac{q^{\langle \alpha^k, \alpha^k\rangle}} {q^{\langle \beta^k, \beta^k \rangle}}  \frac{|\mathrm{R}_\mu^{ss}(\alpha^k, \mathbb{F}_q)|}{|\mathrm{GL}(\alpha^k, \mathbb{F}_q)|}\bigg) X^{\alpha_*},
\end{equation}
where the sum runs over all $r$-tuples of dimension vectors  $\alpha_*=(\alpha^1,\alpha^2,\ldots, \alpha^r)$ with $r\ge 1$, $\alpha^k\in\Delta_\mu^+$ for $1\le k \le r$ and $\alpha^r \ne 0$, 
and $\beta^k=\sum_{i\ge k} \alpha^i$ for $1\le k \le r$.

\begin{dfn} \label{def_refine_kac_function}
The refined Kac functions $A_\mu^{ss}(\lambda_*,q) \in \mathbb{Q}(q)$ for $\lambda_*\in\mathcal{P}^n$ are defined by the following equation:
\begin{equation}\label{def refined kac log}
(q-1)\mathrm{Log}(Q_\Gamma(q, X_{ik})_{1\le i \le n, k\ge 1}) = \sum_{\lambda_*\in\mathcal{P}^n} A_\mu^{ss}(\lambda_*,q) X^{\lambda^{\mathrm{v}}_*},
\end{equation}
where $\mathrm{Log}$ denotes the plethystic logarithm map defined in \eqref{def_Log}.
\end{dfn}

An equivalent definition for $A_\mu^{ss}(\lambda_*,q)$ is given below:
\begin{equation}\label{def_refined_kac_exp}
Q_\mu(q, X_{ik})_{1\le i \leq n, k\geq 1} =  \mathrm{Exp} \bigg(\frac{1}{q-1}\sum_{\lambda_*\in\mathcal{P}^n} A_\mu^{ss}(\lambda_*,q) X^{\lambda^{\mathrm{v}}_*}\bigg).
\end{equation}

For $\alpha=(\alpha_1,\ldots,\alpha_n)\in\mathbb{N}^n$, let $\Lambda_\alpha$ denote the set defined by
\[
	\Lambda_\alpha := \{ (\lambda^1, \ldots, \lambda^n) \in \mathcal{P}^n : |\lambda^i| = \alpha_i  \textrm{ for } 1\le i \le n \}.
\]

With the notation above, substituting $X_{ik}$ with $X_i$ for $1 \leq i \leq n,\,\ k \geq 1$ in~\eqref{def_refined_kac_exp} yields the following identity:
\begin{equation}
A_\mu^{ss}(\alpha,q) = \sum_{\lambda_*\in \Lambda_\alpha} A_\mu^{ss}(\lambda_*,q).
\end{equation}

We propose the following conjecture based on our computations:

\begin{cjc}
Let $\Gamma$ be a connected quiver of infinite representation type; that is, $\Gamma$ is not a Dynkin quiver. Assume that the stability function $\theta$ is non-zero. For any $\lambda_* \in \mathcal{P}^n$, we have:
\begin{enumerate}
    \item $A_\mu^{ss}(\lambda_*,q)$ is a polynomial in $q$; that is, there exists a polynomial $f[t] \in \mathbb{Q}[t]$ such that $A_\mu^{ss}(\lambda_*,q) = f(q)$ for all prime power $q$.
    \item All coefficients of the polynomial $A_\mu^{ss}(\lambda_*,q)$ are non-negative integers.
\end{enumerate}
\end{cjc}

Examples of $A_\mu^{\mathrm{ss}}(\lambda_*, q)$ are provided in the appendix. If $\Gamma$ is a Dynkin quiver, then $A_\mu^{\mathrm{ss}}(\lambda_*, q)$ may not be a polynomial; see Example 1 in the appendix for a counterexample.

The closed formula \eqref{A=sum(H)}, identities \eqref{P=Exp(A/(q-1))} and \eqref{def_refined_kac_exp} may find applications in several areas, including motivic Donaldson–Thomas theory, the study of cohomological Hall algebras,  
the representation theory of quantum groups, and the geometry of Nakajima quiver varieties. 
They also suggest further investigation into wall-crossing phenomena and scattering diagram structures associated with stability conditions. 
The coefficients of the usual Kac polynomials can be interpreted as the dimensions of the homogeneous components of BPS Lie algebras, by Theorem~A of Botta and Davison \cite{B-D 2023}. 
It is therefore natural to expect the existence of a semistable analogue of these BPS Lie algebras, for which a similar statement should hold.

\vspace{0.2cm}
\textbf{\large{Acknowledgments}}
\vspace{0.1cm}

The author is grateful to Xueqing Chen and Bangming Deng for their insightful discussions and valuable feedback on an earlier draft of this paper.

\vspace{3mm}
\textbf{\large{Appendix. Examples}}
\vspace{2mm}

\textbf{Example 1.} 
Let $\Gamma$ be the quiver $A_2$ whose companion matrix is:
\[
\left[
\begin{array}{cc}
0 & 1 \\ 
0 & 0 \\ 
\end{array}
\right],
\]
and $\theta$ be the stability function defined by $\theta((m,n)) = m-n$. Consider the slope $\mu=0$.
Based on identities~\eqref{num_of_semistables}, \eqref{P=Exp(A/(q-1))}, and \eqref{def_refined_kac_exp}, our computations show that:

\begin{center} \label{table_1_1}
	\def\arraystretch{1.2}
	\begin{tabular}{l|l} % <-- Alignments: 1st column left, 2nd middle and 3rd right, with vertical lines in between
		% \textbf{} & \textbf{} & \textbf{}\\
		\hline
		$\ \alpha$ & $A_\mu^{ss}(\alpha, q)$ \\
		\hline
		$(1,1)$ & $1$ \\
		$(2,2)$ & $0$ \\
		$(3,3)$ & $0$ \\
		$(4,4)$ & $0$ \\
		$\cdots$ & $\cdots$ \\
		$\cdots$ & $\cdots$ \\
		$\cdots$ & $\cdots$ \\
		\hline
	\end{tabular}
      \quad \quad \quad
	\def\arraystretch{1.2}
	\begin{tabular}{l|l} % <-- Alignments: 1st column left, 2nd middle and 3rd right, with vertical lines in between
		% \textbf{} & \textbf{} & \textbf{}\\
		\hline
		$\ \lambda_*$ & $A_\mu^{ss}(\lambda_*, q)$ \\
		\hline
		$([1],[1])$ & $1$ \\
		$([1,1],[1,1])$ & $-\frac{1}{q}$ \\
		$([2],[2])$ & $\frac{1}{q}$ \\
		$([1,1,1],[1,1,1])$ & $\frac{1}{q^3}$ \\
		$([2,1],[2,1])$ & $-\frac{1}{q^3} - \frac{1}{q^2}$ \\
		$([3],[3])$ & $\frac{1}{q^2}$ \\
		$\cdots$ & $\cdots$ \\
		\hline
	\end{tabular}
\end{center}

\vspace{0.5cm}
\textbf{Example 2.} 
Let $\Gamma$ be the quiver with companion matrix:
\[
\left[
\begin{array}{cc}
0 & 3 \\ 
0 & 0 \\ 
\end{array}
\right],
\]
and let $\theta$ be the stability function defined by $\theta((m,n)) = m-n$. Consider the slope $\mu=0$.
Our computations show that:

\begin{center} \label{table_2_1}
	\def\arraystretch{1.2}
	\begin{tabular}{l|l} % <-- Alignments: 1st column left, 2nd middle and 3rd right, with vertical lines in between
		% \textbf{} & \textbf{} & \textbf{}\\
		\hline
		$\ \alpha$ & $A_\mu^{ss}(\alpha, q)$ \\
		\hline
		$(1,1)$ & $q^2+q+1$ \\
		$(2,2)$ & $q^5+q^4+3q^3+3q^2+3q+1$ \\
		$(3,3)$ & $q^{10}+q^9+3q^8+5q^7+8q^6+12q^5+17q^4+16q^3+13q^2+6q+2$ \\
		$(4,4)$ & $q^{17}+q^{16}+3q^{15}+5q^{14}+10q^{13}+ 14q^{12}+25q^{11}+35q^{10}+55q^9$\\
		           & $+\,71q^8+97q^7+109q^6+117q^5 +96q^4+68q^3+33q^2+13q+3$ \\
		$\cdots$ & $\cdots$ \\
		\hline
	\end{tabular}
\end{center}
\vspace{0.5cm}
\begin{center} \label{table_2_2}
	\def\arraystretch{1.2}
	\begin{tabular}{l|l} % <-- Alignments: 1st column left, 2nd middle and 3rd right, with vertical lines in between
		% \textbf{} & \textbf{} & \textbf{}\\
		\hline
		$\ \lambda_*$ & $A_\mu^{ss}(\lambda_*, q)$ \\
		\hline
		$([1],[1])$           & $q^2+q+1$ \\
		$([1,1],[1,1])$     & $q^5+q^4+2q^3+2q^2+2q+1$ \\
		$([2],[2])$           & $q^3+q^2+q$ \\
		$([1,1,1],[1,1,1])$        & $q^{10}+q^9+3q^8+5q^7+7q^6+9q^5+11q^4+10q^3$ \\
                                                  & $+\, 9q^2+5q+2$ \\
		$([2,1],[2,1])$             & $q^6+3q^5+5q^4+5q^3+3q^2+q$ \\
		$([3],[3])$                   & $q^4+q^3+q^2$ \\
		$([1,1,1,1],[1,1,1,1])$  & $q^{17}+q^{16}+3q^{15}+5q^{14}+10q^{13}+14q^{12}+24q^{11}$ \\
                                                  & $+\, 32q^{10}+46q^9+55q^8+68q^7+71q^6+72q^5+60q^4$ \\
                                                  & $+\, 45q^3+25q^2+11q+3$ \\
		$([2,1,1],[2,1,1])$        & $q^{11}+3q^{10}+8q^9+15q^8+25q^7+32q^6+34q^5 $ \\
                                                  & $+\, 27q^4+17q^3+7q^2+2q$ \\
		$([2,2],[2,2])$             & $q^9+q^8+3q^7+3q^6+5q^5+3q^4+2q^3$ \\
		$([3,1],[3,1])$             & $q^7+3q^6+5q^5+5q^4+3q^3+q^2$ \\
		$([4],[4])$                   & $q^5+q^4+q^3$ \\
		$\cdots$ & $\cdots$ \\
		\hline
	\end{tabular}
\end{center}

\vspace{0.5cm}
\textbf{Example 3.}
Let $\Gamma$ be the quiver with companion matrix:
\[
\left[
\begin{array}{cc}
1 & 1  \\ 
0 & 1  \\ 
\end{array}
\right],
\]
and let $\theta$ be the stability function defined by $\theta((m,n)) = m-n$. Consider the slope $\mu=0$.
Our computations show that:

\begin{center} \label{table_2_1}
	\def\arraystretch{1.2}
	\begin{tabular}{l|l} % <-- Alignments: 1st column left, 2nd middle and 3rd right, with vertical lines in between
		\hline
		$\ \alpha$ & $A_\mu^{ss}(\alpha, q)$ \\
		\hline
		$(1,1)$ & $q^2$ \\
		$(2,2)$ & $q^5+q^4+q^3$ \\
		$(3,3)$ & $q^{10}+q^9+3q^8+3q^7+3q^6+q^5+q^4$ \\
		$(4,4)$ & $q^{17}+q^{16}+3q^{15}+5q^{14}+8q^{13}+10q^{12}+13q^{11}+12q^{10}+11q^9$ \\
                        & $+\, 7q^8+4q^7+q^6+q^5$ \\
		$\cdots$ & $\cdots$ \\
		\hline
	\end{tabular}
\end{center}
\vspace{0.5cm}
\begin{center} \label{table_2_2}
	\def\arraystretch{1.2}
	\begin{tabular}{l|l} % <-- Alignments: 1st column left, 2nd middle and 3rd right, with vertical lines in between
		% \textbf{} & \textbf{} & \textbf{}\\
		\hline
		$\ \lambda_*$ & $A_\mu^{ss}(\lambda_*, q)$ \\
		\hline
		$([1],[1])$           & $q^2$ \\
		$([1,1],[1,1])$      & $q^5+q^4$ \\
		$([2],[2])$           & $q^3$ \\
		$([1,1,1],[1,1,1])$         & $q^{10}+q^9+3q^8+3q^7+2q^6$ \\
		$([2,1],[2,1])$              & $q^6+q^5$ \\
		$([3],[3])$                   & $q^4$ \\
		$([1,1,1,1],[1,1,1,1])$    & $q^{17}+q^{16}+3q^{15}+5q^{14}+8q^{13}+10q^{12}+12q^{11} $ \\
                                               & $+\,10q^{10}+7q^9+3q^8$ \\
		$([2,1,1],[2,1,1])$         & $q^{11}+2q^{10}+3q^9+3q^8+2q^7 $ \\
		$([2,2],[2,2])$              & $q^9+q^8+q^7$ \\
		$([3,1],[3,1])$              & $q^7+q^6$ \\
		$([4],[4])$                   & $q^5$ \\
		$\cdots$ & $\cdots$ \\
		\hline
	\end{tabular}
\end{center}

\vspace{0.2cm}
%\newpage
\textit{Email address}: \texttt{jiuzhao.hua@gmail.com}

\end{document}